\documentclass[11pt]{amsart}
\pdfoutput=1
\usepackage{amssymb}
\usepackage{amsmath}
\usepackage{amsfonts}

\usepackage{graphicx}
\makeatletter 
 
 \@addtoreset{equation}{section}
\makeatother

\begin{document}

\author{Keiko Kawamuro}
%\address{Department of Mathematics, Rice University, Houston, TX 77005}
%\email{keiko.kawamuro@rice.edu}

\title{On a generalized Jones conjecture}
%\subjclass{?}

\begin{abstract} 
The author withdraws this paper.
\end{abstract}

\maketitle

\end{document}